\documentclass[12pt]{article}
\usepackage{amssymb}
\usepackage{amsfonts}
\usepackage{amsmath}
\usepackage{amsthm}
\usepackage{eucal}
\usepackage{pstricks}
\usepackage{latexsym}

\theoremstyle{plain}

\newtheorem{prop}{Proposition}
\theoremstyle{remark}
\newtheorem*{remark}{Remark}

\numberwithin{theorem}{section}
\numberwithin{prop}{section}
\begin{document}
\title{Combinatorics of Rooted Trees and Hopf Algebras}
\author{Michael E. Hoffman\\
\small Dept. of Mathematics\\[-0.8ex]
\small U. S. Naval Academy, Annapolis, MD 21402\\[-0.8ex]
\small \texttt{meh@usna.edu}}
\date{\small This version:  \today\\
\small MR Classifications:  Primary 05C05,16W30; Secondary 81T15}
\maketitle
\def\de{\delta}
\def\la{\lambda}
\def\si{\sigma}
\def\Si{\Sigma}
\def\De{\Delta}
\def\tilde{\widetilde}
\def\<{\langle}
\def\>{\rangle}
\def\Fix{\operatorname{Fix}}
\def\Orb{\operatorname{Orb}}
\def\Inv{\operatorname{Inv}}
\def\Wex{\operatorname{Wex}}
\def\Ch{\operatorname{Ch}}
\def\rank{\operatorname{rank}}
\def\id{\operatorname{id}}
\def\P{\mathcal P}
\def\A{\mathcal A}
\def\PP{\mathfrak P}
\def\NN{\mathfrak N}
\def\DD{\mathfrak D}
\def\T{\mathcal T}
\def\PT{\mathcal P\mathcal T}
\def\TT{k\{\mathcal T\}}
\def\HK{\mathcal H_K}
\def\HGL{\mathcal H_{GL}}
\begin{abstract}
We begin by considering the graded vector space with a basis 
consisting of rooted trees, with grading given by the count of non-root 
vertices.
We define two linear operators on this vector space, the growth and 
pruning operators,
which respectively raise and lower grading; their commutator is the
operator that multiplies a rooted tree by its number of vertices, and
each operator naturally associates a multiplicity to each pair of rooted
trees.
By using symmetry groups of trees we define an inner product with 
respect to which the growth and pruning operators are adjoint, and obtain
several results about the associated multiplicities.
\par
Now the symmetric algebra on the vector space of rooted trees
(after a degree shift) can be endowed with a coproduct to make a 
Hopf algebra; this was defined by Kreimer in connection with 
renormalization.  We extend the growth and pruning operators, as
well as the inner product mentioned above, to Kreimer's Hopf algebra.
On the other hand, the vector space of rooted trees itself can be given a 
noncommutative multiplication:  with an appropriate coproduct, this leads 
to the Hopf algebra of Grossman and Larson.  We show that the 
inner product on rooted trees leads to an isomorphism of the Grossman-Larson
Hopf algebra with the graded dual of Kreimer's Hopf algebra, correcting an 
earlier result of Panaite.
\end{abstract}
\section{Introduction}
In recent work on renormalization of quantum field theory, D. Kreimer
and his collaborators \cite{BK,CK,K1,K2,K3,Kr} introduce a Hopf algebra 
(here denoted $\HK$) whose generators are rooted trees.  
Various other Hopf algebras based on
rooted trees have appeared in the literature, in particular that of
R. Grossman and R. G. Larson \cite{GL}, which F. Panaite \cite{Pa} 
connected to $\HK$.
The proof of the principal result of \cite{Pa}
actually contains an error due to the confusion of two kinds of 
multiplicities associated to triples of rooted trees.  In this paper
we show how to correct Panaite's result, while clarifying the
combinatorial significance of these multiplicities.
\par
Kreimer's Hopf algebra $\HK$ admits a derivation called the growth operator,
which is important in describing the relation of this algebra to 
another Hopf algebra studied earlier by A. Connes and H. Moscovici
\cite{CM}.
We introduce a complementary derivation called the pruning operator.
In fact, we find it easiest to start (in \S2) in the vector space $k\{\T\}$
of rooted trees rather than in $\HK$.  There we have growth and pruning
operators (denoted $\NN$ and $\PP$ respectively), and for each pair
of rooted trees $t,t'$ with $|t|\le|t'|$ (where $|t|$ is the number of 
vertices of $t$) there are natural multiplicities $n(t;t')$ and $m(t;t')$ 
associated with $\NN$ and $\PP$ respectively.  
A comparison of these multiplicities using symmetry groups of rooted trees 
leads to the definition of an inner product with respect to which $\NN$ 
and $\PP$ are adjoint.
The operators $\NN$ and $\PP$ are very similar to the
adjoint operators that appear in R. Stanley's theory of differential
posets \cite{St2,St3}, and can be described in terms of S. Fomin's 
somewhat more general theory of dual graded graphs
\cite{F0,F1,F2}.  The techniques of Stanley and Fomin can be used
to obtain various results about $n$ and $m$.
\par
In \S3 we extend the growth and pruning operators (as well as the
inner product) to the Hopf algebra $\HK$.  In this setting the growth
and pruning operators are not quite adjoint, but the deviation
from adjointness is easily described (Proposition 3.3).  We also
describe the duals of the growth and pruning operators.  In \S4
we extend the multiplicities $n$ and $m$ to multiplicities 
$n(t_1,t_2;t_3)$ and $m(t_1,t_2;t_3)$ associated to triples of
rooted trees $t_1,t_2,t_3$ with $|t_1|+|t_2|=|t_3|$.  We then give
an explicit isomorphism of the Grossman-Larson Hopf
algebra onto the graded dual of $\HK$ (Proposition
4.4), providing a corrected version of Panaite's result.  We also 
show how the isomorphism gives another description of the duals of 
the growth and pruning operators.
\par
In addition to the references cited above, two recent articles treat
aspects of $\HK$ not covered here:  see 
\cite{B} for connections between Kreimer's Hopf algebra and earlier 
work on Runge-Kutta methods, and \cite{Fo} for an analysis of the 
primitives of $\HK$.  
The author thanks the referee for bringing them to his attention.
\section{The Vector Space of Rooted Trees}
A rooted tree is a partially ordered set (whose elements are called 
vertices) with a unique greatest element (the root vertex), such
that, for any vertex $v$, the vertices exceeding $v$ in the partial
order form a chain.  If $v$ exceeds $w$ in the partial order, we
call $w$ a descendant of $v$ and $v$ an ancestor of $w$.  If $v$
covers $w$ in the partial order (i.e., $v$ is an ancestor of $w$
and there are no vertices between $v$ and $w$ in the order), we call 
$w$ a child of $v$ and $v$ the parent of $w$.
\par
We can visualize a rooted tree as a directed graph by putting an edge from
each vertex to each of its children:  the root is the only vertex with 
no incoming edge.  
We call a vertex terminal if it has no outgoing edges (i.e., no children).  
The condition that the set of ancestors of any vertex forms a chain insures 
this graph has no cycles.
For a finite rooted tree $t$, we denote by $|t|$ the number of vertices 
of $t$: let $\T$ be the set of finite rooted trees, and
$\T_n=\{t\in\T : |t|=n+1\}$.
For example, $\T_0=\{\bullet\}$, where $\bullet$ is the tree consisting
of only the root vertex, and below are the four elements of $\T_3$, with
the root placed at the top:
\vskip .6in
\hskip .5in
\psline{*-*}(4.5,0)(4.5,.5)\psline{*-*}(4.5,.5)(4.5,1)\psline{*-*}(4.5,1)(4.5,1.5)
\psline{*-*}(5.5,.5)(5.75,1)\psline{*-*}(5.75,1)(6,.5)\psline{*-*}(5.75,1)(5.75,1.5)
\psline{*-*}(6.75,1)(7,1.5)\psline{*-*}(7,1.5)(7.25,1)\psline{*-*}(7.25,1)(7.25,.5)
\psline{*-*}(8,1)(8.3,1.5)\psline{*-*}(8.3,1.5)(8.3,1)\psline{*-*}(8.3,1.5)(8.6,1)
\par
We can define a partial order $\preceq$ on the set $\T$ itself by
setting $t\preceq t'$ if $t$ can be obtained from $t'$ by removing some
non-root vertices and edges; of course $t\preceq t'$ implies $|t|\le |t'|$.
Evidently $t'$ covers $t$ for this order exactly in the case when $t$ can be
obtained by removing from $t'$ a single terminal vertex and the edge
into it.  In this case we write $t\lhd t'$.
\par
If $t\lhd t'$, we can get from $t'$ to $t$ by removing a (terminal)
edge, and from $t$ to $t'$ by adding an edge (and accompanying terminal
vertex).  This leads to the definitions of two numbers associated
with the pair $(t,t')$:
$$
n(t;t')=\text{the number of vertices of $t$ to which a new edge
can be added to obtain $t'$,} 
$$
and
$$
m(t;t')=\text{the number of edges of $t'$ which when removed leave $t$} .
$$
That these two numbers are not always equal can be seen from the example
of
\vskip .2in
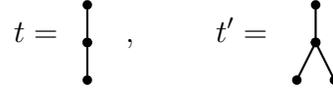
$$
t=\ \psline{*-*}(.25,0)(.25,.5)\psline{*-*}(.25,0)(.25,-.5)
\hskip .3in , \hskip .4in
t'=\ \psline{*-*}(.25,-.5)(.5,0)\psline{*-*}(.5,0)(.5,.5)
\psline{*-*}(.5,0)(.75,-.5)
$$
\vskip .3in
\par\noindent
where $n(t;t')=1$ and $m(t;t')=2$.
\par
The relation between the numbers $n(t;t')$ and $m(t;t')$ can be clarified
by introducing symmetry groups of trees.  For a rooted tree $t$, let
$V(t)$ be its set of vertices:  then for each $v\in V(t)$, there is a
rooted tree $t_v$ consisting of $v$ and its descendants with the
order inherited from $t$.  We call this the subtree of $t$ with $v$ as 
root.  For $v\in V(t)$,
let $SG(t,v)$ be the group of permutations of identical branches out of
$v$, i.e., if $\{v_1,v_2,\dots,v_k\}$ are the children of $v$, then
$SG(t,v)$ is the group generated by the permutations that exchange
$t_{v_i}$ with $t_{v_j}$ when they are isomorphic rooted trees.  The 
symmetry group of $t$ is the direct product
$$
SG(t)=\prod_{v\in V(t)} SG(t,v) .
$$
For a given $v\in V(t)$, let $\Fix (v,t)\le SG(t)$ be the subgroup of
$SG(t)$ that fixes $v$; note that $\Fix (w,t)\le\Fix(v,t)$ whenever $w$
is a descendant of $v$.
\par
Now suppose $t\lhd t'$, and let $v\in V(t)$ be such that, when a new
edge and terminal vertex $w$ are added to $t$ at $v$, the result is $t'$.
If $\Orb(v,t)$ is the orbit of $v$ under $SG(t)$, then evidently
$$
n(t;t')=|\Orb(v,t)|=|SG(t)/\Fix(v,t)|=\frac{|SG(t)|}{|\Fix(v,t)|} .
$$
On the other hand, if $\Orb(w,t')$ is the orbit of $w\in V(t')$ under
$SG(t')$, then 
$$
m(t;t')=|\Orb(w,t')|=|SG(t')/\Fix(w,t')|=\frac{|SG(t')|}{|\Fix(w,t')|} .
$$
But there is an evident identification of $\Fix(v,t)$ with $\Fix(w,t')$,
so we have the following result.
\begin{prop} If $t\lhd t'$, then $|SG(t)|m(t;t')=n(t;t')|SG(t')|$.
\end{prop}
\par
Let 
$$
\TT=\bigoplus_{n\ge 0} k\{\T_n\}
$$ 
be the graded vector space (over a field $k$ of characteristic 0) 
with basis consisting of rooted trees:  we put the rooted tree $t$ 
in grade $|t|-1$.
We define two linear operators on $\TT$ as follows.  
For $n\ge 0$, the growth operator $\NN: k\{\T_n\}\to k\{\T_{n+1}\}$ 
is defined by
\begin{equation}
\NN(t)=\sum_{t\lhd t'} n(t;t')t' ,
\end{equation}
and for $n\ge 1$ the pruning operator $\PP: k\{\T_n\}\to k\{\T_{n-1}\}$ 
is given by
\begin{equation}
\PP(t)=\sum_{t'\lhd t} m(t';t)t' ;
\end{equation}
we set $\PP(\bullet)=0$.  Then $\PP$ and $\NN$ satisfy the following
commutation relation.
\begin{prop} As operators on $\TT$, $[\PP,\NN]=\DD$, where $\DD$ is
the operator given by $\DD(t)=|t|t$.
\end{prop}
\begin{proof} It suffices to show $\PP\NN(t)-\NN\PP(t)=|t|t$ for any
rooted tree $t$.  Let $V(t)=\{v_1,\dots, v_n,v_{n+1},\dots,v_{|t|}\}$
be the vertices of $t$, with $v_i$ terminal for $1\le i\le n$.  Then
$$
\NN(t)=\sum_{i=1}^{|t|} t_i\quad\text{and}\quad
\PP(t)=\sum_{i=1}^n t^{(i)} ,
$$
where $t_i$ is the tree obtained from $t$ by adding a new edge and
terminal vertex to $t$ at $v_i$, and $t^{(i)}$ comes from $t$ by
removing the edge that ends in $v_i$.  Then
$$
\PP\NN(t)=\sum_{i=1}^{|t|}\PP(t_i)=
\sum_{i=1}^{|t|}\left(t+\sum_{1\le j\le n,j\ne i} (t_i)^{(j)}\right)
=|t|t+\sum_{i=1}^{|t|}\sum_{1\le j\le n,j\ne i} (t_i)^{(j)}
$$
and
$$
\NN\PP(t)=\sum_{j=1}^n \NN(t^{(j)})=
\sum_{i=1}^{|t|}\sum_{1\le j\le n,j\ne i} (t^{(j)})_i .
$$
Since $(t_i)^{(j)}=(t^{(j)})_i$ for $i\ne j$, the conclusion follows.
\end{proof}
\par
Now we can endow $\TT$ with an inner product by setting
$$
(t,t')=|SG(t)|\de_{t,t'}
$$
for any rooted trees $t,t'$.
\begin{prop} The operators $\NN$ and $\PP$ are adjoint with respect
to the inner product $(\cdot,\cdot)$.
\end{prop}
\begin{proof} It suffices to show
$$
(\NN(t),t')=(t,\PP(t'))
$$
when $t\lhd t'$ (otherwise both sides are zero).  In this case,
we have
$$
(\NN(t),t')=n(t;t')(t',t')=n(t;t')|SG(t')|
$$
from equation (1) and 
$$
(t,\PP(t'))=m(t;t')(t,t)=m(t;t')|SG(t)|
$$
from equation (2); but then the result follows by Proposition 2.1.
\end{proof}
Putting the last two results together gives the following.
\begin{prop} For rooted trees $t_1$ and $t_2$,
$$
(\NN(t_1),\NN(t_2))-(\PP(t_1),\PP(t_2))=\begin{cases}
0,&\text{if $t_1\ne t_2$,}\\
|t||SG(t)|,&\text{if $t_1=t_2=t.$}\end{cases}
$$
\end{prop}
\begin{proof} This follows from the calculation
$$
(\NN(t_1),\NN(t_2))-(\PP(t_1),\PP(t_2))=(t_1,(\PP\NN-\NN\PP)(t_2))
=(t_1,\DD(t_2))=|t_2|(t_1,t_2) .
$$
\end{proof}
\begin{remark} The second alternative of this result can be written
$$
\sum_{t\lhd t'}n(t;t')^2|SG(t')|-\sum_{t''\lhd t}m(t'';t)^2|SG(t'')|
=|t||SG(t)|,
$$
or, dividing by $|SG(t)|$,
$$
\sum_{t\lhd t'}n(t;t')m(t;t')-\sum_{t''\lhd t}n(t'';t)m(t'';t)=|t|
$$
for any rooted tree $t$.
\end{remark}
We can extend the definitions of $m(t;t')$ and $n(t;t')$ to any
pair of rooted trees $t,t'$ with $|t'|-|t|=k\ge 0$ by setting
\begin{equation}
\NN^k(t)=\sum_{|t'|=|t|+k}n(t;t')t'
\end{equation}
and
\begin{equation}
\PP^k(t')=\sum_{|t'|=|t|+k}m(t;t')t .
\end{equation}
With these definitions, we have the following result.
\begin{prop} Let $t,t'$ be rooted trees with $|t|\le |t'|$.  Then
\newline
1. $n(t;t')|SG(t')|=|SG(t)|m(t;t')$.
\newline
2. If $|t|\le k\le |t'|$, then
$$
n(t;t')=\sum_{|t''|=k} n(t;t'')n(t'';t'),
$$
and similarly for $n$ replaced by $m$.
\newline
3. $n(t;t')$ and $m(t;t')$ are nonzero if and only if $t\preceq t'$.
\end{prop}
\begin{proof}
The first part follows immediately from equations (3) and (4):
$$
n(t;t')|SG(t')|=(t',\NN^{|t'|-|t|}(t))=(\PP^{|t'|-|t|}(t'),t)=m(t;t')|SG(t)| .
$$
For the second part, we have for $|t|\le k\le |t'|$,
\begin{align*}
n(t;t')&=\frac{(\NN^{|t'|-|t|}(t),t')}{|SG(t')|}\\
&=\frac{(\NN^{k-|t|}(t),\PP^{|t'|-k}(t'))}{|SG(t')|}\\
&=\sum_{|t''|=k}\frac{(\NN^{k-|t|}(t),m(t'';t')t'')}{|SG(t')|}\\
&=\sum_{|t''|=k}\frac{(\NN^{k-|t|}(t),t'')}{|SG(t'')|}
\frac{|SG(t'')|}{|SG(t')|}m(t'';t')\\
&=\sum_{|t''|=k}n(t;t'')n(t'';t') .
\end{align*}
(For the corresponding equation with $m$ replacing $n$, reverse the roles
of $\NN$ and $\PP$.)  Finally, the third part is evident for $|t'|-|t|=1$
and can be proved by induction on $|t'|-|t|$ using the second part.
\end{proof}
\par\noindent
\begin{remark} The second and third parts say that $\T$ is a 
weighted-relation poset,
in the terminology of \cite{Ho}, for either of the weights $n(t;t')$ or
$m(t;t')$.  In fact, $\T$ with weights $n(t;t')$ is discussed as 
Example 7 in \cite{Ho}.  In the terminology of \cite{F1}, $\T$ with
multiplicities $n(t;t')$ and $\T$ with multiplicities $m(t;t')$ are
a pair of graded graphs that are $\mathbf r$-dual for the sequence
$\mathbf r=(0,1,2,\dots)$.
\end{remark}
If $t\preceq t'$, we can think of $n(t;t')$ as counting the ways of
building up $t'$ from $t$ by adding new edges and terminal vertices,
and $m(t;t')$ as counting ways of getting from $t'$ to $t$ by removing 
terminal edges.  In particular, since $\bullet\preceq t$ for every 
rooted tree $t$, we can think of $n(\bullet;t)$ as the number of ways
to build up $t$, and $m(\bullet;t)$ as the number of ways to tear it
down.  A more precise formulation can be given using the idea of
labellings of trees:  a labelling of a rooted tree $t$ is a bijection
$f: V(t)\to\{0,1,\dots,|t|\}$ such that $f(v)>f(w)$ whenever $v$ is a 
descendant of $w$ (necessarily $f$ sends the root vertex 
to 0).  We call labellings $f$ and $g$ equivalent if $f\phi=g$ for 
some $\phi\in SG(t)$. 
\begin{prop} Let $t$ be a rooted tree.  Then $t$ has $m(\bullet;t)$ labellings
and $n(\bullet;t)$ labellings mod equivalence.
\end{prop}
\begin{proof}
First note that $|SG(t)|n(\bullet;t)=m(\bullet;t)$ by the first part of
Proposition 2.5 since $\bullet$ has trivial symmetry group.  It follows
from the discussion of \cite[Ex. 7]{Ho} that $n(\bullet;t)$ counts
labellings mod equivalence, and the statement about $m(\bullet;t)$ follows 
since each equivalence class of labellings has $|SG(t)|$ elements. 
\end{proof}
\begin{remark} The ``Connes-Moscovici weight'' \cite{BK,K3} or ``tree
multiplicity'' \cite{B} of $t$ is $n(\bullet; t)$.
Cf. \cite[Sect. 22]{St1} and \cite[Ex. 5.1.4-20]{Kn}, where a hook-length 
formula for the number of labellings of $t$ is given:  this is 
$m(\bullet;t)$.
\end{remark}
\par
In \cite{St3} Stanley defined the notion of a sequentially differential 
poset (generalizing his definition of a differential poset in \cite{St2}).
A sequentially differential poset $P$ is a locally finite graded poset 
with a single element $\hat 0$ in grade 0, so that the linear operators
$$
U(p)=\sum_{p\lhd p'} p'\quad\text{and}\quad D(p)=\sum_{p'\lhd p} p' 
$$
on $k\{P\}$ satisfy the identity $(DU-UD)(p)=r_j p$ for any $p\in P$
of rank $j$:  here $r_0,r_1,\dots$ are nonnegative integers.  The
results of \cite{St3} can be applied to $\T$ (with $r_j=j+1$), provided 
we replace $U$ and $D$ with $\NN$ and $\PP$ respectively, and suitably 
reinterpret the statements of theorems to incorporate multiplicities.
For example, for $x\in P$ Stanley writes $e(x)$ for the number of saturated
chains from $\hat 0$ to $x$, but in the proofs $e(x)$ really appears as
the inner product of $x$ with $U^k\hat 0$, where $k$ is the rank of $x$:
so for a rooted tree $x\in\T_k$ we replace $e(x)$ by
$$
(\NN^k\bullet,x)=n(\bullet;x)(x,x)=n(\bullet;x)|SG(x)|=m(\bullet;x).
$$
Let $w=w_1w_2\cdots w_r$ be a word in $\NN$ and $\PP$, and let $x\in\T_k$.  
Clearly $(w\bullet,x)=0$ unless (a) for each $1\le i\le r$, the number
of $\PP$'s in $w_iw_{i+1}\cdots w_r$ does not exceed the number of $\NN$'s;
and (b) the number of $\NN$'s minus the number of $\PP$'s in $w$ is $k$.
In this case we call $w$ a valid $x$-word, and we have the following result.
\begin{prop} Let $x\in\T_k$, $w=w_1\cdots w_r$ a valid $x$-word.  
Let $S=\{i: w_i=\PP\}$.  For each $i\in S$, let $a_i=|\{j: j\ge i, w_j=\PP\}|$,
$b_i=|\{j: j > i, w_j=\NN\}|$, and $c_i=b_i-a_i$.  Then
$$
(w\bullet,x)=m(\bullet;x)\prod_{i\in S}\binom{c_i+2}{2} .
$$
\end{prop}
\begin{proof} Replace $U,D$ by $\NN,\PP$ in Theorem 2.3 of \cite{St3}.
\end{proof}
This result has the following corollary (cf. Theorem 1.5.2 of \cite{F1}).
\begin{prop} For any rooted tree $x\in\T_k$ and nonnegative integer $a$,
$$
\sum_{|t|=k+a+1}m(x;t)n(\bullet;t)=n(\bullet;x)\prod_{i=2}^{a+1}\binom{k+i}{2}.
$$
\end{prop}
\begin{proof} In Proposition 2.7 set $w=\PP^a\NN^{a+k}$ to get
$$
(\PP^a\NN^{a+k}\bullet,x)=m(\bullet;x)\prod_{i=1}^a\binom{k+a-i+2}{2}
=m(\bullet;x)\prod_{i=2}^{a+1}\binom{k+i}{2}.
$$
Now the left-hand side can be expanded as
\begin{multline*}
(\NN^{a+k}\bullet,\NN^a(x))=\sum_{|t|=k+a+1}n(x;t)(\NN^{a+k}\bullet,t)
=\sum_{|t|=k+a+1}n(x;t)m(\bullet;t)\\
=\sum_{|t|=k+a+1}m(x;t)|SG(x)|n(\bullet;t) ,
\end{multline*}
where we have the first part of  Proposition 2.5 in the last step.  Hence
$$
\sum_{|t|=k+a+1}m(x;t)|SG(x)|n(\bullet;t) =
m(\bullet;x)\prod_{i=2}^{a+1}\binom{k+i}{2},
$$
and dividing by $|SG(x)|$ gives the conclusion.
\end{proof}
\par\noindent
\begin{remark} In the case $x=\bullet$, this result becomes
$$
\sum_{|t|=a+1}m(\bullet;t)n(\bullet;t)=
\sum_{|t|=a+1}n(\bullet;t)^2|SG(t)|=\prod_{i=2}^{a+1}\binom{i}{2} .
$$
Cf. Corollary 1.5.4 of \cite{F1}.
\end{remark}
In \cite[Ex. 7]{Ho} it is shown that $\sum_{|t|=k+1}n(\bullet;t)=k!.$
Further sum formulas involving $n(\bullet;t)$ appear in \cite[Sect. 5]{K3} 
and \cite[Sect. 5]{B}.
A result of \cite{St3} gives a formula for $\sum_{|t|=k+1}m(\bullet;t)$.
To state it we will need some definitions.  Let $\Inv(k)$ be the
set of involutions in the group $\Si_k$ of permutations of $\{1,2,\dots,k\}$.
For $\si\in\Si_k$, call
$i$ a weak excedance of $\si$ if $\si(i)\ge i$; let $\Wex(\si)$ be the set of
weak excedances of $\si$.  For $\si\in\Si_k$ and 
$i\in\{1,\dots,k\}$, let $\eta(\si,i)$ be the number of integers $j$
such that $j<i$ and $\si(j)<\si(i)$.  Then we have the following result.
\begin{prop} With the definitions above,
$$
\sum_{|t|=k+1}m(\bullet;t)=\sum_{\si\in\Inv(k)}\prod_{i\in\Wex(\si)} 
(\eta(\si,i)+1) .
$$
\end{prop}
\begin{proof}
In the proof of Theorem 2.1 of \cite{St3}, replace $U,D$ with $\NN,\PP$:
in the conclusion, this replaces
$\alpha(0\to k)=\sum_{\rank x=k}e(x)$ with $\sum_{|t|=k+1}m(\bullet;t)$.
\end{proof}
For example, a sum over the four involutions $123$, $213$, $132$, and
$321$ of $\Si_3$ gives
$$
\sum_{|t|=4}m(\bullet;t)=1\cdot 2\cdot 3+1\cdot 3+1\cdot 2+1\cdot 1=12.
$$
\section{Kreimer's Hopf Algebra}
In this section we discuss the Hopf algebra $\HK$ defined by D. Kreimer 
and his collaborators \cite{BK,CK,K1,K2,K3,Kr} in connection with 
renormalization.
As an algebra $\HK$ is generated by rooted trees; so as a vector space
$\HK$ is generated by monomials in rooted trees, i.e.,``forests'' of rooted 
trees.
For a rooted tree $t$, we give the corresponding generator degree $|t|$
in $\HK$; in degree 0, $\HK$ is generated by the unit element 1.  
For example, the degree-3 part of $\HK$ is generated as a vector
space by the four elements
\vskip .2in
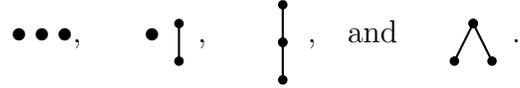
$$
\bullet\bullet\bullet ,
\hskip .3in
\bullet\psline{*-*}(.25,-.25)(.25,.25) 
\hskip .2in , \hskip .3in
\psline{*-*}(.25,-.5)(.25,0)\psline{*-*}(.25,0)(.25,.5)
\hskip .2in \text{,\quad and} \hskip .2in
\psline{*-*}(.25,-.25)(.5,.25)\psline{*-*}(.5,.25)(.75,-.25)
\hskip .4in .
$$
\vskip .2in
\par
There is a linear map $B_+:\HK\to\TT$ which takes a
forest to a single tree with a new root vertex connected to all
the roots of the forest:  e.g.,
\vskip .2in
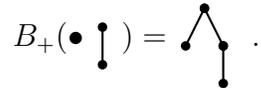
$$
B_+(\bullet\psline{*-*}(.25,-.25)(.25,.25)
\hskip .2in ) =
\psline{*-*}(.25,0)(.5,.5)
\psline{*-*}(.5,.5)(.75,0)
\psline{*-*}(.75,0)(.75,-.5)
\hskip .4in .
$$
\vskip .2in
\par\noindent
The map $B_+$ takes the degree-$n$ part of $\HK$ onto $k\{\T_n\}$:
if we set $B_+(1)=\bullet$, then $B_+$ is a vector space isomorphism.
We write $B_-$ for the inverse of $B_+$.
On the other hand, except for the degree shift, $\HK$ is just the symmetric 
algebra on $\TT$.  
Thus, if $T_n=\dim k\{\T_n\}=|\T_n|$, we have
\begin{equation}
\sum_{n\ge 0}T_nx^n=\prod_{n\ge 1}\frac1{(1-x^n)^{T_{n-1}}}
\end{equation}
from which we can compute recursively $T_0=1$, $T_1=1$, $T_2=2$, $T_3=4$,
$T_4=9$, etc. (see \cite{Sl} for more information).
\par
To define the  bialgebra structure on $\HK$, we let the counit
send all elements of positive degree to 0, and the unit element 1
in degree 0 to $1\in k$.  The comultiplication $\De$ has $\De(1)=1\otimes 1$,
\begin{equation}
\De(t)=t\otimes 1+(\id\otimes B_+)\De(B_-(t))
\end{equation}
for a rooted tree $t$, and $\De(t_1t_2\cdots t_n)=
\De(t_1)\De(t_2)\cdots\De(t_n)$ for monomials $t_1t_2\cdots t_n$.
\par
Equation (6) gives a recursive definition of the coproduct, but there is also
a nonrecursive definition in terms of cuts.  A cut of a rooted tree 
$t$ is a set of edges of $t$.  A cut is elementary if its cardinality is 1.
When the elements of a cut $C$ of $t$ are removed, what remains is a collection
of rooted trees:  the one containing the root is denoted $R^C(t)$, and 
the remaining rooted trees form a monomial denoted $P^C(t)$.  For example,
if $t$ is the tree
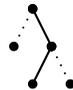
\[
\psline[linestyle=dotted]{*-*}(0,0)(.25,.5)
\psline{*-*}(.25,.5)(.5,0)
\psline{*-*}(.5,0)(.25,-.5)
\psline[linestyle=dotted]{*-*}(.5,0)(.75,-.5)
\]
\vskip .2in
\par\noindent
and $C$ consists of the dotted edges, then 
\vskip .1in
\[
R^C(t)=\hskip .1in\psline{*-*}(0,.5)(0,0)\psline{*-*}(0,0)(0,-.5)
\hskip .2in \text{and}\quad P^C(t)=\bullet\bullet .
\]
\vskip .2in
\par\noindent
The order
of a cut $C$ of $t$ is the largest number of edges in $C$ between the
root of $t$ and any of its terminal vertices:  a cut of order at most 1
is called admissible.  The empty cut $\emptyset$ is the only cut of
order 0 (note that $R^\emptyset(t)=t$ and $P^\emptyset(t)=1$).  The following
formula for $\De(t)$ is proved in \cite{CK}.
\begin{prop} For a rooted tree $t$, $\De(t)$ can be written
$$
\De(t)=t\otimes 1+\sum_{\text{$C$ admissible cut of $t$}}P^C(t)\otimes R^C(t) .
$$
\end{prop}
\par
We can define growth and pruning operators $N$ and $P$ on $\HK$ as
follows.  The growth operator $N$ is simply $\NN$ extended as a derivation,
i.e.,
$$
N(t_1t_2\cdots t_n)=\sum_{i=1}^n t_1\cdots\NN(t_i)\cdots t_n .
$$
We also define $P$ as a derivation, but set $P(t)=\PP(t)$ only
for $|t|\ge 2$; we put $P(\bullet)=1$.  If $D:\HK\to\HK$ is the
extension of $\DD$ as a derivation (i.e., the
linear map that multiplies a monomial by its degree),
then the identity
\begin{equation}
[P,N]=D
\end{equation}
holds.  To prove equation (7), note that both sides are derivations, so 
it suffices to prove it for rooted trees $t$; but in that case (7) follows 
from Proposition 2.2.
The map $B_+$ interacts with the growth and pruning operators as follows.
\begin{prop} For monomials $u$ of $\HK$,
\newline
1. $B_+P(u)=\PP B_+(u)$,
\newline
2. $B_+N(u)=\NN B_+(u)-B_+(\bullet u)$.
\end{prop}
\begin{proof}
Suppose $u=t_1\cdots t_k$ with each $|t_i|\ge 1$.  Then applying
$B_+$ to 
$$
P(u)=P(t_1)t_2\cdots t_k+t_1P(t_2)t_3\cdots t_k+\dots+t_1\cdots t_{n-1}P(t_n)
$$
gives a sum of rooted trees that includes all those obtained by removing
terminal edges of $B_+(u)$, and the cases with $t_i=\bullet$ (hence
$P(t_i)=1$) work correctly:  these are exactly those cases where an
edge coming out of the root of $B_+(u)$ is terminal.  If $u=1$, then
$\PP B_+(1)=\PP(\bullet)=0=B_+P(1)$.  So in any case $B_+P(u)$ coincides
with $\PP B_+(u)$.
\par
Now for $u=t_1\cdots t_k$, $B_+$ applied to
$$
N(u)=N(t_1)t_2\cdots t_k+t_1N(t_2)t_3\cdots t_k+\dots+t_1\cdots t_{n-1}N(t_n)
$$
will include all those trees obtained by adding a new edge to each vertex
of $B_+(u)$ except one--the ``new'' root vertex.  Thus, $B_+N(u)$ is 
missing the term obtained by adding a new edge to the root of $B_+(u)$,
namely $B_+(\bullet u)$.  On the other hand, if $u=1$ we have $B_+N(1)=0
=\NN(\bullet)-B_+(\bullet)$.
\end{proof}
\par
We can extend the inner product of the previous section to $\HK$
by setting
$$
(u_1,u_2)=(B_+(u_1),B_+(u_2))
$$
for monomials $u_1,u_2$; there is no ambiguity since $(B_+(t),B_+(t'))
=(t,t')$ for rooted trees $t,t'$.  With this definition, we can
state the adjointness relation between $P$ and $N$.
\begin{prop} On $\HK$, the adjoint of $P$ with respect to the 
inner product above is $N+M_\bullet$, where $M_\bullet$ is the 
operator that sends $u$ to $\bullet u$; equivalently, the adjoint
of $N$ is $P-\frac{\partial}{\partial\bullet}$.
\end{prop}
\begin{proof} Let $u_1,u_2$ be monomials of $\HK$.  Then 
$$
(u_1,P(u_2))=(B_+(u_1),B_+P(u_2))
=(B_+(u_1),\PP B_+(u_2))
=(\NN B_+(u_1),B_+(u_2)),
$$
from which the first statement follows using the second part of 
Proposition 3.2.
For the second statement,
note that $\frac{\partial}{\partial\bullet}$ is adjoint to $M_\bullet$.
\end{proof}
We now compute the characteristic polynomial of the restriction 
$PN_k$ of $PN$ to the degree-$k$ part of $\HK$.
Let $\Ch(L,\la)=\det(\la I-L)$ for a linear transformation $L$.
\begin{prop}
For $k\ge 1$,
$$
\Ch(PN_k,\la)=\left(\la-\binom{k+1}{2}\right)
\prod_{r=0}^{k-1}\left(\la-\sum_{j=0}^r(k-j)\right)^{T_{k-r}-T_{k-r-1}},
$$
where as above $T_i$ is the dimension of the degree-$i$ part of $\HK$.
\end{prop}
\begin{proof} We follow the argument of \cite[Theorem 4.1]{St2}.  Evidently
$PN_1$ is the identity, so the result holds for $k=1$; assume it inductively
for $k\ge 1$.  From elementary linear algebra
$$
\Ch(NP_{k+1},\la)=\la^{T_{k+1}-T_k}\Ch(PN_k,\la) ,
$$
while from equation (7) we have
$$
\Ch(PN_{k+1},\la)=\Ch(NP_{k+1},\la-(k+1))
$$
since $D_{k+1}=(k+1)I$.  The induction step then follows.
\end{proof}
The preceding result implies that $N_k$ is injective for all $k\ge 1$,
and that $P_k$ is surjective for $k\ge 2$; of course $P_1$ is also surjective.
In addition, the maximal eigenvalue of $PN_k$ is $\binom{k+1}{2}$.  In fact,
the element
$$
f_k=N^{k-1}(\bullet)=\sum_{|t|=k}n(\bullet;t)t
$$
is a corresponding eigenvector.  To see this, note that
$$
PN_k(f_k)=P(f_{k+1})=\sum_{|t'|=k}t'\sum_{|t|=k+1}n(\bullet;t)m(t';t)
=\sum_{|t'|=k}n(\bullet;t')\binom{k+1}{2}t'=\binom{k+1}{2}f_k,
$$
where we have used Proposition 2.8.  The $f_k$ are the ``naturally
grown forests'' of \cite{CK} (where they are denoted $\de_k$).
\par
The following result, which describes how $N$ behaves with respect
to the coproduct, is essentially \cite[Prop. 6]{CK}.  We give
the proof since it can be stated concisely and illustrates the use
of Proposition 3.1.
\begin{prop}
$\De N=(N\otimes\id+\id\otimes N+M_\bullet\otimes D)\De$.
\end{prop}
\begin{proof} Since both sides are derivations, it suffices to show that
$$
\De N(t)=(N\otimes\id+\id\otimes N+M_\bullet\otimes D)\De(t)
$$
for any rooted tree $t$.  As in the proof of Proposition 2.2, write 
$N(t)=\sum_i t_i$, where each $t_i$ is the result of adding an edge to $t$.  
Then
\begin{align*}
\De N(t) &=\sum_it_i\otimes 1+\sum_i\sum_{\text{$C_i$ admissible cut of $t_i$}}
P^{C_i}(t_i)\otimes R^{C_i}(t_i)\\
&=N(t)\otimes 1+\sum_i\sum_{\text{$C_i$ admissible cut of $t_i$}}
P^{C_i}(t_i)\otimes R^{C_i}(t_i) .
\end{align*}
Now each cut $C_i$ of $t_i$ either includes the ``new'' edge or it doesn't.
Suppose first that $C_i$ does not include the new edge. 
Then $C_i$ corresponds to a cut $C$ of $t$ and either 
$P^{C_i}(t_i)\otimes R^{C_i}(t_i)$ is a term in $P^C(t)\otimes NR^C(t)$
(if the new edge is in the component of the root) or a term in 
$NP^C(t)\otimes R^C(t)$ (if it isn't).  Together with the leading term
$N(t)\otimes 1$, these give all the terms of 
$(N\otimes\id+\id\otimes N)\De(t)$.
\par
Now suppose that $C_i$ includes the new edge of $t_i$.  If $C$ is the cut
of $t$ given by $C_i$ minus the new edge, then the new edge must have
been attached to a vertex of $R^C(t)$ (by the definition of admissibility),
and so
$$
P^{C_i}(t_i)\otimes R^{C_i}(t_i)=\bullet P^C(t)\otimes R^C(t) .
$$
Since (for each admissible cut $C$ of $t$) there are $|R^C(t)|$ vertices
to which the new edge could be attached, terms of this form contribute
$(M_\bullet\otimes D)\De(t)$.
\end{proof}
\begin{remark} It follows from this result that the $f_k$, $k\ge 1$,
generate a sub-Hopf-algebra of $\HK$.  This Hopf algebra 
is isomorphic to the graded dual of the universal enveloping algebra of 
$\A^1$, the Lie algebra of formal vector fields on $\mathbf R$ that 
vanish to order 2 at the origin (see \cite{CK}).
\end{remark}
Since $\HK$ is a locally finite commutative Hopf algebra, its graded
dual $\HK^{gr}$ is a locally finite cocommutative Hopf algebra,
hence (by the results of \cite{MM}) the universal enveloping algebra of 
the Lie algebra $\P(\HK^{gr})$,
the primitives of $\HK^{gr}$.  Primitives of $\HK^{gr}$ are dual to
indecomposables of $\HK$, and so are linear combinations of elements $Z_t$
for rooted trees $t$, where $\<Z_t,u\>=\de_{t,u}$ for monomials $u\in\HK$.
The duals of $N$ and $P$ can be described as follows.
\begin{prop} 1. $N^*$ is given by $N^*(Z_\bullet)=0$,
\begin{equation}
N^*(Z_t)=\sum_{|t'|=|t|-1} n(t';t)Z_{t'}
\end{equation}
for $|t|\ge 2$, and 
\begin{equation}
N^*(wv)=(N^*w)v+w(N^*v)+\frac{\partial w}{\partial Z_\bullet}|v|v
\end{equation}
for $w,v\in\HK^{gr}$.
\newline
2. $P^*(w)=Z_\bullet w$ for $w\in\HK^{gr}$.
\end{prop}
\begin{proof}
To prove the statements about $N^*(Z_t)$, note that $\<N^*(Z_t),u\>
=\<Z_t,N(u)\>$ is zero
unless $u$ is a scalar multiple of $t'$, for some $t'\lhd t$; but then
equation (8) follows from equation (1).  Equation (9) follows
from Proposition 3.5 since the multiplication in $\HK^{gr}$ is induced
by $\De$.
\par
For the second part, let $t$ be a rooted tree.  If we write
$P(t)=\sum_i t^{(i)}$ as in the proof of Proposition 2.2, then evidently
$$
\bullet\otimes P(t)=\sum_i \bullet\otimes t^{(i)}
$$
are (by Proposition 3.1) exactly those terms of $\De(t)$ of the form 
$\bullet\otimes t'$.
Now let $u=t_1t_2\cdots t_n$ be a monomial of $\HK$.  
Then
\begin{align*}
\De(u)&=\prod_{i=1}^n\De(t_i)\\
&=\prod_{i=1}^n(1\otimes t_i+\bullet\otimes P(t_i)+\cdots)\\
&=1\otimes t_1t_2\cdots t_n + \bullet\otimes(P(t_1)t_2\cdots t_n+\dots+
t_1\cdots t_{n-1}P(t_n))+\cdots\\
&=1\otimes u+\bullet\otimes P(u)+\cdots
\end{align*}
and thus
$$
\<Z_\bullet w,u\>=\<Z_\bullet\otimes w, \De(u)\>=\<w,P(u)\>=\<P^*(w),u\>
$$
for all $w\in\HK^{gr}$ and monomials $u$ of $\HK$.
\end{proof}
\begin{remark} The Lie algebra $\P(\HK^{gr})$ is in fact free; see
\cite{Fo}.
\end{remark}
\section{The Grossman-Larson Hopf Algebra}
We can define a noncommutative multiplication on the graded vector space
$\TT$ as follows.  Let $t,t'$ be rooted trees, and suppose $B_-(t)=t_1t_2
\cdots t_k$.  There are $|t'|^k$ rooted trees obtainable by attaching each
of the $k$ rooted trees $t_1,t_2,\dots,t_k$ to some vertex of $t'$ (by
a new edge): 
let $t\circ t'\in\TT$ be the sum of these trees (If $t=\bullet$, we define
$t\circ t'$ to be $t'$).  For example,
\vskip .2in
$$
\psline{*-*}(.25,0)(.5,.5)\psline{*-*}(.5,.5)(.75,0)
\hskip .4in \circ
\psline{*-*}(.25,0)(.25,.5)
\hskip .3in = \hskip .1in
\psline{*-*}(.25,0)(.5,.5)\psline{*-*}(.5,.5)(.5,0)\psline{*-*}(.5,.5)(.75,0)
\hskip .4in + \quad 2
\psline{*-*}(.25,0)(.5,.5)\psline{*-*}(.5,.5)(.75,0)\psline{*-*}(.75,0)(.75,-.5)
\hskip .5in +
\psline{*-*}(.25,-.5)(.5,0)\psline{*-*}(.5,0)(.5,.5)\psline{*-*}(.5,0)(.75,-.5)
$$
\vskip .2in
\par\noindent
while
\vskip .2in
$$
\psline{*-*}(.25,0)(.25,.5)
\hskip .3in \circ
\psline{*-*}(.25,0)(.5,.5)\psline{*-*}(.5,.5)(.75,0)
\hskip .4in = \hskip .1in
\psline{*-*}(.25,0)(.5,.5)\psline{*-*}(.5,.5)(.5,0)\psline{*-*}(.5,.5)(.75,0)
\hskip .4in + \quad 2
\psline{*-*}(.25,0)(.5,.5)\psline{*-*}(.5,.5)(.75,0)\psline{*-*}(.75,0)(.75,-.5)
\hskip .4in .
$$
\vskip .2in
This product makes $\TT$ a graded algebra:  note that for
$t\in k\{\T_n\}$ and $t'\in k\{\T_m\}$, we have $t\circ t'\in k\{\T_{n+m}\}$.
The element $\bullet\in\T_0$ is a two-sided identity.
Note also that
\vskip .1in
$$
B_+(\bullet)\circ t=
\psline{*-*}(.25,.25)(.25,-.25)
\hskip .2in\circ\ t=\NN(t)
$$
\vskip .2in
\par\noindent
for any rooted tree $t$.
\par
Now define a coproduct $\De:\TT\to\TT\otimes\TT$ by 
\begin{equation}
\De(t)=\sum_{I\cup J=\{1,\dots,k\}}B_+(t_I)\otimes B_+(t_J)
\end{equation}
where $B_-(t)=t_1\cdots t_k$ and the sum is over pairs $(I,J)$ of
(possibly empty) subsets $I,J$ of $\{1,\dots,k\}$ such that 
$I\cup J=\{1,\dots,k\}$:  $t_I$ means the product of $t_i$ for
$i\in I$.  The following result is proved in \cite{GL} and
\cite{GVF}:  the main things to check are the associativity of
the product $\circ$ (Lemma 2.6 of \cite{GL}) and the compatibility
of the coproduct with $\circ$ (Lemma 2.8 of \cite{GL}).
\begin{prop} The vector space $\TT$ with product $\circ$ and coproduct
$\De$ is a graded Hopf algebra $\HGL$.
\end{prop}
Since the coproduct $\De$ is cocommutative, by results of \cite{MM} it 
follows that $\HGL$ is the universal enveloping algebra
on its Lie algebra $\P(\HGL)$ of primitives.  
From equation (10), elements of the form $B_+(t)$, where $t$ is a rooted 
tree, are primitive.
We call such elements ``primitive trees'':  they are those rooted trees 
whose root has exactly one child.
If we let $\PT$ be the set of primitive trees (graded, like $\T$, by
the number of non-root vertices), then we have following result 
(for another proof, see \cite[Theorem 4.1]{GL}).
\begin{prop} The vector space $k\{\PT\}$ generated by the primitive
trees is $\P(\HGL)$.
\end{prop}
\begin{proof} Since
$$
B_+(t_1)\circ B_+(t_2)= B_+(t_1t_2) + B_+(B_+(t_1)\circ t_2),
$$
$k\{\PT\}\subseteq\P(\HGL)$ is a sub-Lie-algebra.
Also, since $B_+$ is an isomorphism of $k\{\T_{n-1}\}$ onto
$k\{\PT_n\}$, we have $\dim k\{\PT_n\}=T_{n-1}$.
Then the Poincar\'e-Birkhoff-Witt theorem implies that the universal 
enveloping algebra of $k\{\PT\}$ has the same dimension in grade 
$n$ as does the symmetric algebra on $k\{\PT\}$:  but in view of equation 
(5), this is $T_n=\dim (\HGL)_n$.  Hence $k\{\PT\}=\P(\HGL)$.
\end{proof}
\par
Suppose $t_1,t_2,t_3$ are rooted trees so that $|t_1|+|t_2|=|t_3|$.
If there is an elementary cut $C$ of $t_3$ so that
\begin{equation}
P^C(t_3)=t_1\quad\text{and}\quad R^C(t_3)=t_2 ,
\end{equation}
let $m(t_1,t_2;t_3)$ be the number of distinct elementary cuts $C$
of $t_3$ for which equations (11) hold:  otherwise, set $m(t_1,t_2;t_3)=0$.
If (and only if) $m(t_1,t_2;t_3)\ne 0$, it is
also true that $t_3$ can be obtained by attaching (via a new edge)
the root vertex of $t_1$ to some vertex of $t_2$:  let $n(t_1,t_2;t_3)$
be the number of vertices of $t_2$ for which this is true.  Evidently
$$
n(\bullet,t_2;t_3)=n(t_2;t_3)\quad\text{and}\quad
m(\bullet,t_2;t_3)=m(t_2;t_3)
$$
for trees $t_2\lhd t_3$, so we have generalized the multiplicities of \S2.
(The reader is warned that $n(t_1,t_2;t_3)$ as used in \cite{CK} and 
\cite{Fo} is our $m(t_1,t_2;t_3)$.)
We now show how symmetry groups can be used to relate the two multiplicities.
\begin{prop} For rooted trees $t_1,t_2,t_3$ with $|t_1|+|t_2|=|t_3|$,
$$
|SG(t_1)||SG(t_2)|m(t_1,t_2;t_3)=n(t_1,t_2;t_3)|SG(t_3)|.
$$
\end{prop} 
\begin{proof}
For any rooted tree $t$ and $v\in V(t)$, let $\Fix(t_v,t)
\le SG(t)$ be the subgroup of $SG(t)$ that holds $t_v$ (the subtree
of $t$ with $v$ as root) pointwise fixed.  We can assume there is an
elementary cut $C=\{e\}$ of $t_3$ so that equations (11) hold (otherwise
both sides of the conclusion are zero).  If $e$ has source $v$ and
target $w$, then $t_w$ is isomorphic to $t_1$.  Also, if $\Orb (e,t_3)$
is the orbit of $e$ under $SG(t_3)$, then
$$
m(t_1,t_2;t_3)=|\Orb(e,t_3)|=|SG(t_3)/\Fix(t_v,t_3)\times SG(t_w)|=
\frac{|SG(t_3)|}{|\Fix(t_v,t_3)||SG(t_1)|}.
$$
On the other hand, since $R^C(t_3)$ is isomorphic to $t_2$,
$$
n(t_1,t_2;t_3)=|\Orb(v,R^C(t_3))|=|SG(R^C(t_3))/\Fix(v,R^C(t_3))|=
\frac{|SG(t_2)|}{|\Fix(v,R^C(t_3))|}.
$$
Since there is an evident identification of $\Fix(t_v,t_3)$ with
$\Fix(v,R^C(t_3))$, we have
$$
\frac{|SG(t_1)|}{|SG(t_3)|}m(t_1,t_2;t_3)=\frac{n(t_1,t_2;t_3)}{|SG(t_2)|}
$$
and the conclusion follows.
\end{proof}
We can now use the inner product on Kreimer's Hopf algebra $\HK$ 
to define an isomorphism of $\HGL$ onto the graded dual of $\HK$.
\begin{prop}
There is an isomorphism $\chi:\HGL\to\HK^{gr}$ defined by
$$
\<\chi(t),u\>=(B_-(t),u)=(t,B_+(u))
$$
for any rooted tree $t$ and monomial $u$ of $\HK$.
\end{prop}
\begin{proof}
Since $\HK$ is locally finite, it suffices to prove that $\chi$
is an injective homomorphism.
We first show $\chi$ is a homomorphism, i.e., that
\begin{multline*}
\<\chi(t_1\circ t_2),u\>=\<\chi(t_1)\otimes\chi(t_2),\De(u)\>=\\
\sum_u\<\chi(t_1),u'\>\<\chi(t_2),u''\>=
\sum_u (t_1,B_+(u'))(t_2,B_+(u'')) 
\end{multline*}
for any monomial $u$ of $\HK$ with coproduct
\begin{equation}
\De(u)=\sum_u u'\otimes u'' .
\end{equation}
In view of Proposition 4.2, $\HGL$ is generated as an algebra by
the primitive trees.  So it suffices to show that
\begin{equation}
\<\chi(B_+(t)\circ t_2),u\>=\sum_u (B_+(t),B_+(u'))(t_2,B_+(u''))
=\sum_u (t,u')(t_2,B_+(u'')).
\end{equation}
Now from the definition of $n(t_1,t_2;t_3)$, 
\begin{multline*}
\<\chi(B_+(t)\circ t_2),u\>=(B_+(t)\circ t_2,B_+(u))=
\sum_{|t_3|=|t|+|t_2|} n(t,t_2;t_3)(t_3,B_+(u)) =\\
n(t,t_2;B_+(u))|SG(B_+(u))| .
\end{multline*}
On the other hand, if $\De(u)$ is given by equation (12), then
\begin{equation}
\De(B_+(u))=B_+(u)\otimes 1+\sum_u u'\otimes B_+(u'')
\end{equation}
by equation (6).  Now the only nonzero terms of 
$$
\sum_u (t,u')(t_2,B_+(u''))
$$
are those with $u'=t$ and $t_2=B_+(u'')$:  and (comparing Proposition 3.1
with equation (14)) there are $m(t,t_2;B_+(u))$ such terms.  
Hence
$$
\sum_u (t,u')(t_2,B_+(u''))=m(t,t_2;B_+(u))|SG(t)||SG(t_2)|,
$$
and equation (13) follows from Proposition 4.3:  thus, $\chi$ is a 
homomorphism.
\par
Now suppose $v=\sum_i a_it_i\in\ker\chi$.  Then 
$$
\<\chi(v),u\>=\sum_ia_i(t_i,B_+(u))=0
$$
for all monomials $u$ of $\HK$.  But setting $u=B_-(t_i)$ implies that
$a_i=0$ for each $i$, so $v=0$.
\end{proof}
\par\noindent
\begin{remark}
In \cite[Prop. 2.1]{Pa} (and also in \cite[Theorem 14.16]{GVF}) it is 
wrongly asserted that the map sending $B_+(t)$ to $Z_t$ induces an 
isomorphism of $\HGL$ onto $\HK^{gr}$:  the error is due to a failure
to distinguish the multiplicities $n(t_1,t_2;t_3)$ and $m(t_1,t_2;t_3)$,
since Panaite confuses the coefficients in 
$$
[Z_{t_1},Z_{t_2}]=\sum_{|t_3|=|t_1|+|t_2|} (m(t_1,t_2;t_3)-m(t_2,t_1;t_3))Z_{t_3}
$$
with those in 
$$
B_+(t_1)\circ B_+(t_2)-B_+(t_2)\circ B_+(t_1)=\sum_{|t_3|=|t_1|+|t_2|}
(n(t_1,t_2;t_3)-n(t_2,t_1;t_3))B_+(t_3) .
$$
In fact, since
$$
\<\chi(B_+(t)),u\>=(t,u)=|SG(t)|\de_{t,u}=|SG(t)|\<Z_t,u\> ,
$$
we have $\chi(B_+(t))=|SG(t)|Z_t$.
\end{remark}
We can use the isomorphism $\chi$ to express the duals of $P$ and
$N$ as maps of $\HGL$ (cf. Proposition 3.6 above).
\begin{prop} 1. The map $\chi^{-1}P^*\chi:\HGL\to\HGL$ is left
multiplication by $B_+(\bullet)$, i.e., $\chi^{-1}P^*\chi(t)=\NN(t)=
B_+(\bullet)\circ t$.
\newline
2. For rooted trees $t$, $\chi^{-1}N^*\chi(t)=\PP(t)-B_+\frac{\partial}
{\partial\bullet}B_-(t)$.
\end{prop}
\begin{proof}
Using Propositions 2.3, 3.2, and 3.3, we have for any rooted tree $t$ and 
monomial $u$ of $\HK$, 
$$
\<\chi(t),P(u)\>=(t,B_+P(u))=(t,\PP B_+(u))=(\NN(t),B_+(u))=\<\chi(\NN(t)),u\>
$$
and
\begin{multline*}
\<\chi(t),N(u)\>=(t,B_+N(u))=(t,\NN B_+(u))-(t,B_+(\bullet u))=
(\PP(t),B_+(u))-(B_-(t),\bullet u)\\
=(\PP(t),B_+(u))-\left(\frac{\partial}{\partial\bullet}B_-(t),u\right)
=\left\<\chi\left(\PP(t)-B_+\frac{\partial}{\partial\bullet}B_-(t)\right),u\right\> .
\end{multline*}
\end{proof}

\end{document}